\documentclass{amsart}

\newtheorem{thm}{Theorem}[section]
\newtheorem{lem}[thm]{Lemma}

\theoremstyle{definition}

\newtheorem{conj}[thm]{Conjecture}

\theoremstyle{remark}

\numberwithin{equation}{section}

\begin{document}

\title{The 2-log-convexity of the Ap\'{e}ry Numbers}

\author{William Y. C. Chen}
\address{Center for Combinatorics, LPMC-TJKLC
Nankai University
 Tianjin 300071, P. R. China}
\email{chen@nankai.edu.cn}
\thanks{ The authors  wish to thank  the referee, Tomislav  Do\v{s}li\'{c} and
  Tanguy Rivoal
 for helpful comments.
 This work was supported by  the 973
Project, the PCSIRT Project of the
 Ministry of Education,  and the National Science
Foundation of China.}

\author{Ernest X. W. Xia}
\address{Center for Combinatorics, LPMC-TJKLC
Nankai University
 Tianjin 300071, P. R. China}
\email{xxw@cfc.nankai.edu.cn}


\date{September 4, 2009 and, in revised form, May 11, 2010.}



\begin{abstract}
 We present
 an approach to proving
 the 2-log-convexity of
 sequences satisfying
     three-term recurrence
   relations. We show
    that the Ap\'{e}ry numbers, the
    Cohen-Rhin numbers,   the
     Motzkin numbers, the Fine numbers,
      the Franel numbers of order
 $3$ and $4$ and the
 large Schr\"{o}der numbers
  are all 2-log-convex.
Numerical evidence suggests that all these sequences are
  $k$-log-convex  for any $k\geq 1$  possibly except for
  a constant number of terms at the beginning.
\end{abstract}

\maketitle

\section{Introduction}

\allowdisplaybreaks

In his proof of the irrationality of $\zeta(2)$
 and $\zeta(3)$,
Ap\'{e}ry  \cite{Apery} introduced the
 following
  numbers $A_n$ and $B_n$ as given by
\begin{align}
 A_n&=\sum_{k=0}^n
 {n \choose k}^2
 {n+k\choose k}^2,
 \label{A}\\[6pt]
B_n&=\sum_{k=0}^n {n \choose k}^2
 {n+k\choose k}.\label{B}
\end{align}
The numbers $A_n$  and $B_n$
are often called the Ap\'{e}ry numbers.
It has been shown by Ap\'{e}ry  \cite{Apery}  that
 $A_n$ and $B_n$
  satisfy the following
   three-term recurrence relations
   for $n\geq 2$,
\begin{align}
A_n&=\frac{34n^3-51n^2+27n-5}{n^3}A_{n-1}
-\frac{(n-1)^3}{n^3}A_{n-2},
 \label{R-1}\\[6pt]
B_n&=\frac{11n^2-11n+3}{n^2}B_{n-1}
 +\frac{(n-1)^2}{n^2}B_{n-2},
\label{R-2}
\end{align}
where $A_0=1,\, A_1=5, \, B_0=1,\, B_1=3$;
 see also  \cite{Koepf,Str}.
 Congruences of the Ap\'ery numbers have been
 investigated by Ahlgren, Ekhad, Ono,
 and Zeilberger \cite{A},
  Beukers \cite{Beukers1985,Beukers1987}, Chowla  and
  Clowes\cite{Chowla}
 and Gessel \cite{Gessel}.
  Note that the recurrence relations \eqref{R-1} and
 \eqref{R-2}  can be
 derived  by using Zeilberger's algorithm \cite{Wilf}.

 Cohen  \cite{Cohen} and Rhin
 obtained  the following
 recurrence relation of the numbers $U_n$
    in connection   with the rational approximation of  $\zeta(4)$,
     see also \cite{K},
\begin{align}\label{U}
U_{n+1}= R(n)U_n+ G(n)U_{n-1}, \qquad n\geq 1,
\end{align}
where $U_0=1$, $U_{1}=12$ and
\[
R(n)=\frac{3(2n+1) (3n^2+3n+1)(15n^2+15n+4)} {(n+1)^5},\quad
G(n)=\frac{3n^3(3n-1)(3n+1)}{(n+1)^5}.
\]
Expressions of $U_n$ as double sums of products of binomial
coefficients
 have been derived by
Krattenthaler and Rivoal \cite{K}  and
 Zudilin \cite{Zudilin,W-Zudilin}.

In this paper, we shall establish the 2-log-convexity of the
sequences of the Ap\'ery numbers $A_n$, $B_n$,
 the Cohen-Rhin numbers $U_n$ and some
  other combinatorial sequences based on the
   three-term recurrence relations. Recall that
  an infinite positive sequence
  $\{a_n\}_{n=0}^\infty$ is said to be
   log-convex if for all $n\geq 1$,
   \begin{equation} \label{an}
   a_n^2\leq a_{n-1}a_{n+1}.
   \end{equation}
   We say that
    $\{a_n\}_{n=0}^\infty$ is 2-log-convex
     if $\{a_n\}_{n=0}^\infty$ is log-convex
      and for all $n\geq 1$,
  \begin{align}\label{2}
   \left(a_{n}a_{n+2}-a_{n+1}^2\right)^2
    \leq \left(a_{n-1}a_{n+1}-a_{n}^2\right)
   \left(a_{n+1}a_{n+3}-a_{n+2}^2\right).
  \end{align}
   Meanwhile, the
   sequence $\{a_n\}_{n=0}^\infty$
   is called strictly log-convex (2-log-convex)
   if the inequality in \eqref{an} (\eqref{2})
    is strict
   for all $n\geq 1$.
   Do\v{s}li\'{c} \cite{D} proved
   the log-convexity of $A_n$ by induction. In fact, using
    similar arguments  one can
     show that $\{B_n\}_{n=0}^\infty$ and
      $\{U_n\}_{n=0}^\infty$  are log-convex.

This paper is organized as follows.
 In Section 2, we give a general
  framework to prove the
 2-log-convexity of a
 sequence $\{S_n\}_{n=0}^\infty$  based on
     a lower bound $f_n$ and
     an upper bound $g_n$ for the
     ratio $S_n/S_{n-1}$,
      where the numbers $S_n$ satisfy
      a three-term
       recurrence relation.
      Section 3 demonstrates how to find the bounds $f_n$ and $g_n$.
        Section 4 is devoted to the computations of
        the upper bounds for the ratios $A_n/A_{n-1}$, $B_n/B_{n-1}$
         and $U_n/U_{n-1}$. In Section 5,
          we show that the sequences of $A_n$, $B_n$, $U_n$,
      the Motzkin numbers, the Fine numbers,
      the Franel numbers of order
 $3$ and $4$ and the
 large Schr\"{o}der numbers are all 2-log-convex.
   We conclude this paper with
   a conjecture on the infinite log-convexity in the spririt of the infinite log-concavity introduced
 by Moll \cite{Moll}.

\section{A  criterion}

In this section,
we present a criterion
 for the $2$-log-convexity
 of
a sequence $\{S_n\}_{n=0}^\infty$
 satisfying a three-term recurrence
relation.  We need the assumption
 that the ratio $S_n/S_{n-1}$ has a
lower bound $f_n$ and an upper bound $g_n$.

\begin{thm}\label{L-lemma-1}
Suppose $\{S_n\}_{n=0}^\infty$ is a  positive log-convex
  sequence
 that satisfies  the
recurrence
 relation
  \begin{align}\label{R-S}
  S_n=b(n)S_{n-1} +c(n)S_{n-2}
  \end{align}
for $n\geq 2$. Let
\begin{align*}
a_3(n)=&2b(n+2)b^2(n+1)+2b(n+1)c(n+2) -b^3(n+1)
\\[6pt]
&-b(n+1)b(n+2)b(n+3) -b(n+3)c(n+2) -c(n+3)b(n+1),
\\[6pt]
a_2(n)=&4b(n+1)b(n+2)c(n+1) +2c(n+1)c(n+2)+b^2(n+1)b(n+2)b(n+3)\\[6pt]
&+b(n+1)b(n+3)c(n+2) +b^2(n+1)c(n+3)-3c(n+1)b^2(n+1)
 \\[6pt]
&-b(n+3)b(n+2)c(n+1) -c(n+3)c(n+1)-b^2(n+2)b^2(n+1) \\[6pt]
&-2b(n+2) b(n+1)c(n+2)-c^2(n+2),\\[6pt]
a_1(n)=&-c(n+1)\big(2b(n+2)c(n+2)-2b(n+2)c(n+1)
\\[6pt]
&-2b(n+3)b(n+2)b(n+1)
-b(n+3)c(n+2)-2c(n+3)b(n+1) \\[6pt]
& +3c(n+1)b(n+1)+2b^2(n+2)b(n+1) \big),
\\[6pt]
a_0(n)=&-c^2(n+1) \left(c(n+1)-b(n+2)b(n+3)
 -c(n+3)+b^2(n+2)\right)
\end{align*}
and
\[
\Delta(n)=4a_2^2(n)-12a_1(n)a_3(n).
\]
Assume that
  $a_3(n)<0$ and $\Delta(n)>0$ for all
  $n\geq N$, where $N$ is a positive integer.
   If there exist
 $f_n$ and $g_n$ such that for
  all $n\geq N$,\\[6pt]
($C_1$)
  $f_n\leq \frac{S_n}{S_{n-1}} < g_n$;
\\[6pt]
($C_2$) $f_n\geq \frac{-2a_2(n)-\sqrt{\Delta(n)}}{6a_3(n)}$;
\\[6pt]
($C_3$) $a_3(n)g_n^3+a_2(n)g_n^2+a_1(n)g_n+a_0(n)>0$,
\\[6pt]
 then $\{S_n\}_{n=N}^\infty$ is
   strictly 2-log-convex, that is,
   for $n \geq N$,
  \begin{align}\label{2-logconvex}
   \left(S_{n-1}S_{n+1}-S_n^2\right)
  \left(S_{n+1}S_{n+3}-S_{n+2}^2\right)
   >\left(S_nS_{n+2}-S_{n+1}^2\right)^2.
  \end{align}
\end{thm}

\noindent {\it Proof.} By the recurrence relation \eqref{R-S},
 we have
 \begin{eqnarray*}
\lefteqn{\left(S_{n-1}S_{n+1}-S_n^2\right)
\left(S_{n+1}S_{n+3}-S_{n+2}^2\right)
-\left(S_nS_{n+2}-S_{n+1}^2\right)^2}\\[6pt]
&& =S_{n+1}\left(2S_{n}S_{n+1}S_{n+2} +S_{n-1}S_{n+1}S_{n+3}
-S_{n+1}^3-S_n^2S_{n+3} -S_{n-1}S_{n+2}^2\right)\\[6pt]
&&=S_{n+1}\left(a_3(n)S_{n}^3+a_2(n)S_{n}^2S_{n-1}+a_1(n)S_nS_{n-1}^2
+a_0(n)S_{n-1}^3\right).
 \end{eqnarray*}
Since $\{S_n\}_{n=0}^\infty$
 is a positive sequence, in order to prove
  \eqref{2-logconvex}, it suffices to show that for all $n\geq N$,
  \begin{align}\label{aim}
a_3(n)\left(\frac{S_n}{S_{n-1}}\right)^3
+a_2(n)\left(\frac{S_n}{S_{n-1}}\right)^2
 +a_1(n)\frac{S_n}{S_{n-1}}
+a_0(n)>0.
  \end{align}
Consider the polynomial $f(x)=a_3(n)x^3+a_2(n)x^2 +a_1(n)x+a_0(n)$.
 Note that
 \[
f'(x)=3a_3(n)x^2+2a_2(n)x+a_1(n).
 \]
Since $a_3(n)<0$ and $\Delta(n)>0$ for all $n\geq N$,
 we see that the quadratic
function $f'(x)$ is negative for
$x>\frac{-2a_2(n)-\sqrt{\Delta(n)}}{6a_3(n)}$. Thus,
 $f(x)$ is  strictly decreasing on the interval
 $[\frac{-2a_2(n)
 -\sqrt{\Delta(n)}}{6a_3(n)},+\infty
 )$. From the assumption
  $g_n > f_n\geq \frac{-2a_2(n)
 -\sqrt{\Delta(n)}}{6a_3(n)}$,
  it follows that
  $f(x)$ is  strictly decreasing
  on the interval $[f_n,g_n]$.
  Since $\frac{S_{n}}{S_{n-1}} \in [f_n,g_n]$, it remains to show that
    $f(g_n)>0$ for any $n\geq N$,
    which is equivalent to  condition ($C_3$), that is,
    \[
a_3(n)g^3_n + a_2(n)g^2_n+a_1(n)g_n+a_0(n)>0
    \]
for any $n\geq N$.
 This completes the proof. \qed

\section{A heuristic approach to
 computing the bounds}

In this section, we present a procedure
 to derive a lower bound $f_n$ and an
upper bound $g_n$ for the ratio
 $S_n/S_{n-1}$ based on a three-term
recurrence relation of $S_n$. We first describe how to obtain an
upper bound $g_n$ as required in Theorem \ref{L-lemma-1}. As will be
seen, this procedure is not guaranteed to give an upper bound $g_n$,
but it is practically valid for many cases.

 Assume that  $\lim\limits_{n\rightarrow \infty}
b(n)=b$ and $\lim\limits_{n\rightarrow \infty} c(n)=c$,
 where $b$ and $c$ are two constants
  and $b^2+4c> 0$.
  All sequences   considered in this paper
  satisfy this condition.
   Let
   \begin{align}
   \label{X}
    x_0=\frac{b+\sqrt{b^2+4c}}{2}.
   \end{align}
        We begin with  the case  $c(n)<0$, and we shall
        try to construct $g_n$ which
            satisfies the
             condition $(C_3)$ together with
           the following inequality:
     \begin{align}\label{condition}
       g_{n+1}-\left(b(n+1)
       +\frac{c(n+1)}{g_n}\right)>0.
     \end{align}
     In fact,
         the condition
       (\ref{condition}) is essential to find
       an upper bound $g_n$ for $S_n/S_{n-1}$. As will be seen in the
       following lemma, if we find a function
        $g_n$ satisfying (\ref{condition})
         and $S_n/S_{n-1}<g_n$ for small $n$,
       then we can deduce that $g_n$ is an
       upper bound for $S_n/S_{n-1}$ for any
        $n$.

\begin{lem} \label{L-1}
Let $S_n$  be the sequence  defined by the recurrence relation
\eqref{R-S}. Assume that $N$ is a positive integer such that
 $c(n)<0$
    for $n\geq N$.
   If $ \frac{S_{N}}{S_{N-1}}
    \leq g_{N} $ and the  condition \eqref{condition}
   holds for
    $n\geq N$,  then we have
   for $n\geq N$,
  \begin{align}\label{case-1}
\frac{S_n}{S_{n-1}}\leq g_n.
  \end{align}
\end{lem}

\noindent {\it Proof.}
 We use induction on $n$.
  Obviously, the lemma holds for
    $n=N$.
 We assume that
  it
   is true for $n=m\geq N$,
   that is,
$ \frac{S_m}{S_{m-1}} <g_m $. Since
  $c(m)<0$
 for $m\geq N$, we see that
\begin{align}\label{M-1}
c(m+1)\frac{S_{m-1}}{S_m}<\frac{c(m+1)}{g_m}.
\end{align}
We now consider the case $n=m+1$.
 From \eqref{R-S}
 and
\eqref{M-1} it follows that
\begin{align}\label{M-2}
\frac{S_{m+1}}{S_m}  =b(m+1)+c(m+1)
 \frac{S_{m-1}}{S_m} \leq
b(m+1)+\frac{c(m+1)}{g_m}.
\end{align}
From  \eqref{condition} and
 \eqref{M-2} we deduce that for $m\geq
N$,
\[
g_{m+1}-\frac{S_{m+1}}{S_m}\geq g_{m+1} - \left(
b(m+1)+\frac{c(m+1)}{g_m}\right)> 0,
\]
which is the statement of the
 lemma  for $n=m+1$. This
 completes the proof. \qed

Now we present a heuristic procedure
to find the  desired upper bound $g_n$.
Let $g_n=x_0$ as given by \eqref{X}.
If $g_n$ satisfies
 the conditions $(C_3)$
 and \eqref{condition}, then
  $g_n$ is the desired choice.
 Otherwise, let $g_n=x_0+\frac{x}{n}$.
 Substitute $g_n$ into \eqref{condition}
 and let $Y(n)$  denote the numerator
 of the left hand side of
   \eqref{condition},
   which is often a polynomial in $n$
   and $x$. Setting
  the coefficient
   of the highest degree in $n$
   of $Y(n)$ to be $0$, we obtain
   an equation
    in $x$. If $x_1$ is the unique
     solution of
      this equation, then we set
      $g_n=x_0+\frac{x_1}{n}$.
  If $g_n=x_0+\frac{x_1}{n}$ satisfies
   the conditions $(C_3)$  and
   \eqref{condition},
 then $g_n$ is the desired choice.
 Otherwise,
  set $g_n= x_0+\frac{x_1}{n}
  +\frac{x}{n^2}$
   and repeat the above process to find a solution
    $x_2$ of the equation. By iteration, we may find      $x_0,x_1,\ldots,x_i$
       such that $g_n=x_0+\frac{x_1}{n}
       +\frac{x_2}{n^2}+\cdots
       +\frac{x_i}{n^i}$ satisfies the conditions $(C_3)$ and \eqref{condition}.

         For example,
         let $S_n=A_n$, where $A_n$
         is Ap\'{e}ry number defined
         by \eqref{A}.
         Since $\lim\limits_{n\rightarrow \infty}
b(n)=34$ and $\lim\limits_{n\rightarrow \infty}
 c(n)=-1$, by the definition of $A_n$,
 we  have
 $x_0=17+12\sqrt{2}$. Since
 $g_n=17+12\sqrt{2}$ does not satisfy the condition
 $(C_3)$ in Theorem \ref{L-lemma-1},
  we further consider $g_n=17+12\sqrt{2}+\frac{x}{n}$.
  Let $Y(n)$ denote   the  numerator
  of the left hand side of
   \eqref{condition}.
It is easy to see that $Y(n)$ is a cubic polynomial in $n$ with the
leading coefficient equal to
   \[
   E_1=-(17\sqrt{2}-24)
   (48x+864\sqrt{2}+1224).\]
   Setting $E_1=0$  gives
   $x_1=-\frac{51}{2}-18\sqrt{2}$. Again,
   $g_n=x_0+\frac{x_1}{n}$ does not satisfy
   \eqref{condition}. So we continue to consider
   $g_n=x_0+\frac{x_1}{n}
   +\frac{x}{n^2}$ and we find that
   $x_2=\frac{609}{64}\sqrt{2}+\frac{27}{2}$.
    Now, $g_n=x_0+\frac{x_1}{n}+\frac{x_2}{n^2}$
     does not satisfy the condition $(C_3)$. Repeating the
     above procedure, we find that
 $x_3=-\frac{225}{128}
 \sqrt{2}-\frac{645}{256}$ and
 $g_n=x_0+\frac{x_1}{n}
+\frac{x_2}{n^2}+\frac{x_3}{n^3}$ satisfies
 \eqref{condition} and the condition $(C_3)$.

For the case $c(n)>0$, we aim to
construct an upper
 bound $g_n$ which satisfies
   condition $(C_3)$
  and  the following  inequality
     \begin{align}\label{con}
g_n-\left(b(n)+\frac{c(n)} {b(n-1) +\frac{c(n-1)}{g_{n-2}}
}\right)>0.
     \end{align}
     Similarly, if we find a function
        $g_n$ satisfying (\ref{con}) and
        $S_n/S_{n-1}<g_n$ for certain $n$, then we can deduce
        that $g_n$ is an upper bound  for any $n$. To be precise,
        we have the following lemma.

      \begin{lem} \label{L-2}
       Let $S_n$ be  defined by \eqref{R-S}.
        If   there exists a positive
          integer
          $N$ such that
           the inequality  \eqref{con} holds,
           $\frac{S_{N}}{S_{N-1}} \leq g_{N}$,
 $\frac{S_{N+1}}{S_{N}} \leq g_{N+1}$
 and $c(n)>0 $ for $n\geq N$,
   then
  we have
   for $n\geq N$,
  \begin{align}\label{case-2}
 \frac{S_n}{S_{n-1}}\leq g_n.
  \end{align}
\end{lem}

\noindent {\it Proof.}
 We conduct induction on $n$.
  Clearly, the lemma holds
   for $n=N$ and $n=N+1$.
  Assume that it is true
 for  $n=m-2\geq N$, that is,
\begin{align}\label{B-1}
\frac{S_{m-2}}{S_{m-3}} \leq  g_{m-2}.
\end{align}
We shall show that the lemma is true for
 $n=m$, that is,
\begin{align}\label{B-2}
\frac{S_m}{S_{m-1}}\leq  g_{m}.
\end{align}
Since $c(n)>0$ for $n\geq N$, from  \eqref{R-S}
 and \eqref{B-1} it follows
 that
\begin{align}\label{M2-1}
\frac{S_m}{S_{m-1}} & =b(m)+c(m)\frac{S_{m-2}} {S_{m-1}}
=b(m)+\frac{c(m)} {b(m-1)+c(m-1) \frac{S_{m-3}}{S_{m-2}}}
 \\[6pt]
&\leq b(m)+\frac{c(m)}
 {b(m-1)+\frac{c(m-1)}{g_{m-2}}}.
  \nonumber
\end{align}
In view of \eqref{con}
 and \eqref{M2-1}, we find that
 \begin{align*}
g_{m}-\frac{S_m}{S_{m-1}} \geq g_m-\left( b(m)+\frac{c(m)}
 {b(m-1)+\frac{c(m-1)}{g_{m-2}}}\right)
 > 0,
 \end{align*}
 which yields \eqref{B-2}.
 This completes the proof.
 \qed

Now we can use the same approach as
in the case $c(n)<0$
 to find an upper bound $g_n$.
 Moreover, if we have obtain an approximation  $g_n$
 that does not simultaneously satisfy \eqref{condition}
 (\eqref{con}) and the condition $(C_3)$, instead of
 going further to update the estimation of $g_n$,
  we may try to  adjust some coefficients
   to find a desired bound.  For example, let $S_n=B_n$, where
     $B_n$ is defined by
 \eqref{B}. At some point, we get
 \begin{align}\label{G-1}
g_n=&\frac{11}{2} +\frac{5\sqrt{5}}{2}
-\left(\frac{11}{2}+\frac{5\sqrt{5}}{2}
 \right)\frac{1}{n}\\[6pt]
&+\left(\frac{7}{10}\sqrt{5}+\frac{3}{2}\right)\frac{1}{n^2}
 +\frac{1}{25n^3}+\left(\frac{1}{50}+\frac{23\sqrt{5}}{1250}\right)
\frac{1}{n^4} . \nonumber
 \end{align}
Here $g_n$ satisfies the condition $(C_3)$
 in Theorem \ref{L-lemma-1}, but it fails to
 satisfy \eqref{con}.
 If we replace  the coefficient
 $\frac{1}{50}$  in \eqref{G-1} by
 $\frac{1}{25}$,
 then the adjusted bound $g'_n$ satisfies both   conditions $(C_3)$ and
 \eqref{con}.

 To conclude this section, we need to mention that it is much
 easier to find a lower bound   $f_n$ for the ratio $S_n/S_{n-1}$.
    In many cases,
   we have  $f(n)=b(n)$ when $b(n)$ and $c(n)$
    are  positive for $n\geq N$ and   $f_n=b(n)+c(n)$
     when $c(n)$ is negative and $S_{n}\geq
     S_{n-1}$ for $n\geq N$.

\section{Upper bounds for $A_n/A_{n-1}$,
 $B_n/B_{n-1}$
 and $U_n/U_{n-1}$ }

In this section, we shall use the heuristic
approach described in
the previous section to find
 upper bounds for the ratios $A_n/A_{n-1}$,  $B_n/B_{n-1}$
 and $U_n/U_{n-1}$.

\begin{lem}\label{Lem-1}
 Let
\begin{align}\label{P}
P(n) = & 17+12\sqrt{2}
-\left(\frac{51}{2}
+18\sqrt{2}\right)\frac{1}{n}
\\[6pt]
&\qquad  +\left(\frac{27}{2} +\frac{609}{64}
\sqrt{2}\right)\frac{1}{n^2} -\left(\frac{645}{256} +
\frac{225\sqrt{2}}{128}\right)\frac{1}{n^3}.\nonumber
\end{align}
For $n\geq 2$, we have
$ \frac{A_n}{A_{n-1}} <
 P(n).
$
\end{lem}

\noindent {\it Proof.} For the Ap\'{e}ry numbers $A_n$, we use Lemma
\ref{L-1} by setting $N=2$
 and $g_n=P(n)$. Evidently,
  $\frac{A_2}{A_1}<P(2)$.
Also, it is easily checked
 that
\begin{align*}
P(n+1)&
 -\left(\frac{(2n+1)(17n^2+17n+5)} {(n+1)^3}
 -\frac{n^3}{(n+1)^3P(n)}
 \right)\\[6pt]
 &=\frac{9(17-12\sqrt{2})
 (5664n^2-3560\sqrt{2}n+1225)}{
 256(256n^3-384n^2-60\sqrt{2}n
 +288n+90\sqrt{2}-165)(n+1)^3},
\end{align*}
 which is positive
  for $n\geq 2$.
 By lemma \ref{L-1}, we see that $P(n)$ is an upper bound for
 $A_n/A_{n-1}$ when $n\geq 2$. This completes the proof.
 \qed

\begin{lem}\label{Lem-2} Let
\begin{align} \label{T}
T(n)=&\frac{11}{2} +\frac{5\sqrt{5}}{2}
-\left(\frac{11}{2}+\frac{5\sqrt{5}}{2} \right)\frac{1}{n}
 \\[6pt]
&\qquad +\left(\frac{7}{10}\sqrt{5}+\frac{3}{2}\right)\frac{1}{n^2}
 +\frac{1}{25n^3}+\left(\frac{1}{25}+\frac{23\sqrt{5}}{1250}\right)
\frac{1}{n^4}. \nonumber
\end{align}
For $n\geq 20$, we have $ \frac{B_n}{B_{n-1}} < T(n). $
\end{lem}

\noindent {\it Proof.}
 Set   $N=20$ and
  $g_n=T(n)$ in Lemma \ref{L-2}. It is easy to check that
    $\frac{B_{20}}{B_{19}}<T(20)$
    and $\frac{B_{21}}{B_{20}}<T(21)$.
Moreover, it is not difficult to verify that
\begin{align*}
T(n)&-\left(\frac{11n^2-11n+3}{n^2} +\frac{(n-1)^2}
 {n^2\left(
\frac{11n^2-33n+25}{(n-1)^2}
 +\frac{(n-2)^2}{(n-1)^2}
 \frac{1}{T(n-2)}\right)}\right)\\[6pt]
 &\qquad \qquad \qquad =\frac{(123\sqrt{5}-275)J(n)}
 {1250n^4K(n)},
\end{align*}
where $J(n)$ and $K(n)$ are given by
\begin{align*}
J(n)=& 1718750n^6-4656250\sqrt{5}n^5
-18026250n^5+98010000n^4\\[6pt]
&+38885750\sqrt{5}n^4 -136205250\sqrt{5}n^3
-310595950n^3+248642319\sqrt{5}n^2
\\[6pt]
&+557184100n^2-233557457\sqrt{5}n
 -522290000n+199152500+89063225\sqrt{5},
\\[6pt]
K(n)=&2500n^6 -30000n^5+150000n^4 -500\sqrt{5}n^4
-401100n^3+4500\sqrt{5}n^3\\[6pt]
&+642325n^2 -30881\sqrt{5}n^2-619575n
+78143\sqrt{5}n-60525\sqrt{5}+278125.
\end{align*}
It follows  that $J(n)$ and $K(n)$
are positive
 for $n\geq 20$. Hence we have
\begin{align}\label{B2-2}
\frac{11n^2-11n+3}{n^2} +\frac{(n-1)^2} {n^2\left(
\frac{11n^2-33n+25}{(n-1)^2}
 +\frac{(n-2)^2}{(n-1)^2}
 \frac{1}{T(n-2)}\right)} <T(n).
\end{align}
In view of Lemma \ref{L-2}, we deduce that
$T(n)$ is an upper bound for $B_n/B_{n-1}$ when $n\geq 20$. \qed

Using the same procedure, we find the following upper bound for $U_n/U_{n-1}$.
 The proof  is omitted.

\begin{lem}\label{Lem-3} Let
\begin{align}\label{Q}
Q(n) =&135+78\sqrt{3}-\left(\frac{675}{2}
+195\sqrt{3}\right)\frac{1}{n}
+\left(\frac{9737}{48}\sqrt{3}+351\right)\frac{1}{n^2}
\\[6pt]
&\qquad
-\left(\frac{3497}{32}\sqrt{3}+\frac{6045}{32}\right)\frac{1}{n^3}
+\left(\frac{841763}{27648}\sqrt{3}
+\frac{2701}{32}\right)\frac{1}{n^4}.\nonumber
\end{align}
For $n\geq 100$, we have $ \frac{U_n}{U_{n-1}} <Q(n). $
\end{lem}

\section{The 2-log-convexity}

Based on the criterion given
in Theorem \ref{L-lemma-1}
 and the upper bounds obtained in the previous section, we shall give  the
proofs of the 2-log-convexity of the sequences of Ap\'ery numbers
and other aforementioned combinatorial numbers.

\begin{thm}\label{Thm-1}
The sequence $\{A_n\}_{n=0}^\infty$
 is strictly $2$-log-convex.
\end{thm}

\noindent {\it Proof.}  We first consider the case
 $n\geq 2$.
 To apply  Theorem
 \ref{L-lemma-1},
 let
 \[
b(n)=\frac{34n^3-51n^2+27n-5}{n^3}\quad \mbox{and} \quad
c(n)=-\frac{(n-1)^3}{n^3}.
 \]
 It is straightforward to check that $a_3(n)<0$
 and $\Delta(n)>0$
  for $n\geq 2$.     Since
  \[
{n-1 \choose k}^2{n-1+k \choose k}^2\geq {n-2\choose
k}^2{n-2+k\choose k}^2,
  \]
  we have $A_{n-1}\geq A_{n-2}$. Let
\[
f_n=\frac{33n^3-48n^2+24n-4}{n^3}.
\] Thus, by the recurrence relation
  \eqref{R-1},
   we see that
\begin{align}\label{f-1}
\frac{A_n}{A_{n-1}}= & \frac{34n^3-51n^2+27n-5}{n^3}-\frac{(n-1)^3}{n^3}
\frac{A_{n-2}}{A_{n-1}}
\\[6pt]
\geq & \frac{34n^3-51n^2+27n-5-(n-1)^3} {n^3}=
f_n.\nonumber
\end{align}
Set $g_n=P(n)$, where $P(n)$ is given by \eqref{P}. We proceed to verify the
 conditions $(C_1)$, $(C_2)$ and
 $(C_3)$ in  Theorem  \ref{L-lemma-1}.
 By \eqref{f-1} and  Lemma \ref{Lem-1},
  we find that   $f_n \leq \frac{A_{n}}{A_{n-1}}<g_n$,
   which is the condition    $(C_1)$. Define $R_1(n)=6a_3(n)f_n+2a_2(n)$.
   It is easily checked that $R_1(n)=-4\frac{H_1(n)}{L_1(n)}$,
where $H_1(n)$ and $L_1(n)$ are polynomials in $n$
 and the leading  coefficients of $H_1(n)$
  and $L_1(n)$ are positive.
   Hence we   deduce that
   $R_1(n)<0$ for $n\geq 2$. Similarly, define $R_2(n)=\Delta(n)-R_1^2(n)$, which
   can be rewritten as $-96\frac{H_2(n)}{L_2(n)}$
where $H_2(n)$ and $L_2(n)$ are polynomials in $n$
 and the leading coefficients
  of $H_2(n)$
  and $L_2(n)$ are  positive. Consequently, we   deduce
   $R_2(n)<0$ for $n\geq 2$.
  It follows that  for $n\geq 2$,
\[
6a_3(n)f_n+2a_2(n)<-\sqrt{\Delta(n)},
\]
which is equivalent to the following inequality  for $n\geq 2$:
\[
f_n>\frac{-2a_2(n)-\sqrt{\Delta(n)}}{6a_3(n)}.
\]
This is exactly the  condition $(C_2)$. Finally, it remains
 to verify the condition
 $(C_3)$.  To this end, we find that
 \begin{align}\label{I}
a_3(n)g_n^3&+a_2(n)g_n^2+a_1(n)g_n+a_0(n)\\[6pt]
&=9\left(30733178557+21731638968\sqrt{2}\right)
 \frac{H_3(n)}{L_3(n)}, \nonumber
 \end{align}
where $H_3(n)$ and $L_3(n)$ are polynomials in $n$. Observe that the
leading coefficients
of $H_3(n)$
  and $L_3(n)$ are both positive.
  This implies that the right hand side of \eqref{I} is positive
    for $n\geq 2$. Now we are left with the case $n=1$, that is
    \[(A_0A_2-A_1^2)(A_2A_4-A_3^2)
>(A_1A_3-A_2^2)^2,\]
 which can be easily checked.
    This completes the proof.
   \qed

\begin{thm}\label{Thm-2}
The sequence $\{B_n\}_{n=0}^\infty$
 is strictly $2$-log-convex.
\end{thm}

\noindent {\it Proof.}
  For $n\geq 20$,
  apply Theorem
 \ref{L-lemma-1} with
\[
f_n=\frac{11n^2-11n+3}{n^2},
\]
and  $g_n=T(n)$, where $T(n)$ is given
 by \eqref{T}.
 Using the argument in the proof of
   Theorem \ref{Thm-1}, we find that
 $f_n$ and $g_n$ satisfy all the conditions in Theorem
  \ref{L-lemma-1}. Finally, it is easy to verify that for $1\leq n\leq 19$,
     \[\left(
   B_{n-1}B_{n+1}-B_n^2\right)
   \left(
   B_{n+1}B_{n+3}-B_{n+2}^2\right)>\left(
   B_{n}B_{n+2}-B_{n+1}^2\right)^2.
\]
  This completes the proof.
  \qed

\begin{thm}\label{Thm-3}
The sequence $\{U_n\}_{n=0}^\infty$
 is strictly $2$-log-convex.
\end{thm}

The above theorem follows from Theorem \ref{L-lemma-1} by
setting
\[
f_n=\frac{3(2n-1)(3n^2-3n+1)(15n^2-15n+4)}{n^5}
\]
and  setting $g_n=Q(n)$, where $Q(n)$ is given by \eqref{Q}.
  The  proof is similar to that of
 Theorem \ref{Thm-1}, and  it is omitted.

Do\v{s}li\'{c} \cite{D,D2} has proved the log-convexity of several
well-known sequences of combinatorial numbers such as the
  Motzkin numbers $M_n$, the Fine numbers $F_n$, the Franel numbers $F_n^{(3)}$ and $F_n^{(4)}$
  of order
 $3$ and $4$, and the large Schr\"{o}der numbers $s_n$. Based on the
recurrence relations satisfied by these numbers, we utilize Theorem
 \ref{L-lemma-1}
to deduce that these sequences are all strictly 2-log-convex
possibly except for a fixed number of terms at the beginning.

We conclude this paper with a conjecture concerning the infinite
log-convexity of the A\'{e}ry numbers. The notion of infinite
log-convexity is analogous to that of infinite log-concavity
introduced by Moll  \cite{Moll}.
  Given a sequence
$A=\{a_i\}_{0 \leq i \leq \infty}$, define the  operator $\mathcal
{L}$
 by \[ \mathcal{L}(A)
 =\{b_i\}_{0\leq i \leq \infty},\] where
 $
b_i=a_{i-1}a_{i+1}-a_i^2$ for $ i\geq 1$.
 We say that $\{a_i\}_{0 \leq i
 \leq \infty}$ is $k$-log-convex
  if
$\mathcal {L}^j\left(\{a_i\}_{0 \leq i \leq \infty}\right)$ is
 log-convex for $j=0, 1,\ldots, k-1$,
 and that $\{a_i\}_{0 \leq i \leq \infty}$
is $\infty$-log-convex
 if $\mathcal {L}^k\left(\{a_i\}_{0 \leq i
\leq \infty}\right)$ is
 log-convex for any $k \geq 0$.

  \begin{conj}\label{conj}
 The sequences $\{A_n\}_{n=0}^\infty$,
   $\{B_n\}_{n=0}^\infty$,
   $\{U_n\}_{n=0}^\infty$ and $\{s_n\}_{n=0}^\infty$ are infinitely
   log-convex.
   The sequences
    $\{M_n\}_{n=0}^\infty$,
     $\{F_n\}_{n=0}^\infty$,
    $\{F_n^{(3)}\}_{n=0}^\infty$ and
      $\{F_n^{(4)}\}_{n=0}^\infty$ are $k$-log-convex for any $k \geq 1 $
      except for a constant number (depending on $k$)
       of terms at the beginning.
  \end{conj}

\bibliographystyle{amsplain}

\end{document}